\documentclass[a4paper]{article}
\usepackage{amsfonts}

\newtheorem{theorem}{Theorem}[section]
\newtheorem{lemma}[theorem]{Lemma}
\newtheorem{example}[theorem]{Example}
\newtheorem{proposition}[theorem]{Proposition}
\newtheorem{corollary}[theorem]{Corollary}
\newtheorem{definition}[theorem]{Definition}
\newtheorem{remark}[theorem]{Remark}

\newcommand{\LGH}{L^2(G;{\cal H})}
\newcommand{\pil}{\mu_{\lambda,T}}
\newcommand{\B}{{\cal B}}
\newcommand{\id}{{\rm id}}

\begin{document} 
\title{Approximation property of $C^*$-algebraic Bundles}
\author{Ruy Exel\footnote{Partially supported by CNPq, Brazil}\hspace{1.5mm} 
and Chi-Keung Ng}
\date{\today}
\maketitle
\begin{abstract}
In this paper, we will define the reduced cross-sectional 
$C^*$-algebras of $C^*$-algebraic bundles over locally compact groups 
and show that if a $C^*$-algebraic bundle has the approximation 
property (defined similarly as in the discrete case), then the full 
cross-sectional $C^*$-algebra and the reduced one coincide. 
Moreover, if a semi-direct product bundle has the approximation property
and the underlying $C^*$-algebra is nuclear, then the cross-sectional 
$C^*$-algebra is also nuclear. 
We will also compare the approximation property with the 
amenability of Anantharaman-Delaroche in the case of discrete groups.
\end{abstract}
\par\medskip
\par\medskip
\noindent {\small 1991 Mathematics Subject Classifiction: 46L55; 46L05; 46L45}

\par\medskip
\par\medskip
\par\medskip

\noindent {\bf \Large 0\quad Introduction}
\par\medskip
\par\medskip

\noindent $C^*$-algebraic bundles were defined and studied by Fell in [7]. 
In a recent paper ([6]), the first named author studied $C^*$-algebraic 
bundles over discrete groups and defined an interesting notion of 
approximation property. 
He showed that if a $C^*$-algebraic bundle has the approximation property, 
then the full cross-sectional $C^*$-algebra and the reduced one coincide. 
This can be regarded as a kind of amenability for $C^*$-algebraic
bundles over discrete groups. 
\par\medskip

The main objective of this paper is to define and study the approximation
property of $C^*$-algebraic bundles in the case of locally compact 
groups and show that when the bundle is a semi-direct product, this 
gives a good candidate for the notion of amenable group actions. 
We will also give some interesting relations between coactions of 
discrete groups and gradings. 
\par\medskip

In the first section, we will recall some basic materials of $C^*$-algebraic bundles
and give some technical lemmas.
We will also define the canonical coactions 
on the full cross-sectional $C^*$-algebras $C^*(\cal B)$. 
\par\medskip

In section 2, we will define and show the existence of the reduced 
cross-sectional $C^*$-algebra $C^*_r(\cal B)$ (of any $C^*$-algebraic 
bundle $\cal B$). 
In fact, we will give two different ways of defining the reduced 
cross-sectional $C^*$-algebra and show that they coincide. 
We will study the relation between $C^*(\cal B)$ and 
$C^*_r(\cal B)$ as well as the reduced coaction on $C^*_r(\cal B)$. 
\par\medskip

In section 3, we will define approximation property of 
$C^*$-algebraic bundles similar to that of the discrete case in [6]. 
We will then show that if a $C^*$-algebraic bundle $\cal B$ has 
the approximation property, then $C^*({\cal B})=C^*_r(\cal B)$.
\par\medskip

In the final section, we will consider two interesting special cases 
of $C^*$-algebraic bundles. 
The first one is group actions on $C^*$-algebras. 
We show that if the $C^*$-algebra is nuclear and the action has the 
approximation property, then the crossed product is also nuclear
(a generalisation of the case of actions by amenable groups). 
As an application, we show that the approximation property  
is stronger than the amenability of Anantharaman-Delaroche 
in the case of nuclear $C^*$-algebras and discrete groups. 
Moreover, they coincide in some special cases (see Corollary 4.a). 
The second special case is the case of discrete groups. 
We will study the situation when a $C^*$-algebraic bundle over a 
discrete group comes from a group coaction. 
\par\medskip

This article replaces the paper ``A Note on the Representation Theory of
Fell Bundles'' (math. OA/9904013) by the first named author, as well as
a paper entitled ``Approximation properties of $C^*$-algebraic Bundles''
by the second named author, both of which circulated as preprints. 
\par\medskip
\par\medskip
\par\medskip

\section{Preliminary and basic results}
\par\medskip
\par\medskip

We begin this section by recalling the definition and notation 
about $C^*$-algebraic bundles.
First of all, we refer the reader to [7, VIII.2.2] for the notion 
of Banach algebraic bundle $\cal B$ over a locally compact group $G$. 
Let $B$ be the bundle space of $\cal B$ and $\pi$ be the bundle 
projection. 
Let $B_t = \pi^{-1}(t)$ be the fiber over $t\in G$. 
It is clear that $B_e$ is a $C^*$-algebra if $\cal B$ is a $C^*$-algebraic 
bundle (see [7, VIII.16.2] for a definition).
We will use the materials from [7, VIII] implicitly. 
Following the notation of [7], we denote by $\cal L(B)$ the set of all 
continuous cross-sections on $\cal B$ with compact support. 
Moreover, for any $f\in \cal L(B)$, ${\rm supp}(f)$ is the closed support 
of $f$.
Furthermore, let $({\cal L}_p(\mu; {\cal B}), \|\ \|_p)$ be the normed
space as defined in [7, II.15.7] (where $\mu$ is the left Haar 
measure on $G$). 
For simplicity, we will denote ${\cal L}_p(\mu; \cal B)$ by 
${\cal L}_p(\cal B)$. 
By [7, II.15.9], $\cal L(B)$ is dense in ${\cal L}_p(\cal B)$. 
We also need the theory of operator valued integration from [7, II], 
especially, we would like to draw the readers' attention to 
[7, II.5.7] and [7, II.\S 16]. 
\par\medskip

Throughout this paper, $\cal B$ is a $C^*$-algebraic bundle 
over a locally compact group $G$ with bundle space $B$ and bundle
projection $\pi$. 
Denote by $C^*(\cal B)$ the cross-sectional $C^*$-algebra of 
$\cal B$ (see [7, VIII.17.2]).
We recall from [7, VIII.5.8] that there exists a canonical map $m$ 
from the bundle space $B$ to the set of multipliers of ${\cal L}_1(\cal B)$ 
(or $C^*(\cal B)$). 
\par\medskip

\begin{lemma} \label{1.4}
The map $m$ from $B$ to $M(C^*(\cal B))$ is {\it faithful} in the sense that 
if $m_a = m_b$, then either $a=b$ or $a=0_r$ and $b=0_s$ ($r,s\in G$). 
\end{lemma}
\noindent {\bf Proof}: 
Suppose that $\pi(a) = \pi(b) = r$. 
Then $m_{a-b} = 0$ will imply that $a-b = 0_r$ (since $\cal B$ has a strong 
approximate unit and enough continuous cross-sections). 
Suppose that $\pi(a)=r\neq s=\pi(b)$. 
Then there exists a neighbourhood $V$ of $e$ such that 
$rV\cap sV = \emptyset$. 
For any $f\in \cal L(B)$, $a(f(r^{-1}t)) = m_a(f)(t) = m_b(f)(t) = 
b(f(s^{-1}t))$. 
Now let $b_i$ be a strong approximate unit of $\cal B$ and $\{ f_i \}$ be 
elements in $\cal L(B)$ such that ${\rm supp}(f_i)\subseteq V$ and $f_i(e) = b_i$. 
Therefore $ab_i = a(f(e)) = b(f(s^{-1}r)) =0$ and hence $a=0_r$. 
Similarly, $b=0_s$.
\par\medskip

From now on, we will identify $B_r$ ($r\in G$) with its image 
in $M(C^*({\cal B}))$. 
\par\medskip

Let ${\cal B}\times G$ be the $C^*$-algebraic bundle over $G\times G$
with the Cartesian product $B\times G$ as its bundle space such that
the bundle projection $\pi'$ is given by $\pi'(b,t) = (\pi(b),t)$ ($b\in
B$; $t\in G$). 
It is not hard to see that any non-degenerate representation $T'$ of 
$\B\times G$ is of the form $T'_{(b,t)} = T_bu_t$ ($b\in B$; $t\in G$) 
for a non-degenerate representation $T$ of $\B$ and a unitary
representation $u$ of $G$ with commuting ranges.
This gives following lemma.
\par\medskip

\begin{lemma} \label{1.5}
$C^*({\cal B}\times G) = C^*({\cal B})\otimes_{\max}C^*(G)$. 
\end{lemma}

Consider the map $\delta_{\cal B}$ from $B_r$ to 
$M(C^*({\cal B})\otimes_{\max}C^*(G))$ given by $\delta_{\cal B} (b) = 
b\otimes \Delta_r$ where $\Delta_r$ is the canonical image of $r$ in
$M(C^*(G))$. 
Denote again by $\delta_{\cal B}$ the integral form of $\delta_{\cal
B}$. 
Then we have the following equalities.

\begin{eqnarray*}
\delta_{\cal B}(f)(1\otimes k)(g\otimes l)(r,s) & = &
\int_G f(t)g(t^{-1}r)[\int_G k(u)l(u^{-1}t^{-1}s)du]dt\\ 
& = & \int_G\int_G f(t)g(t^{-1}r)k(t^{-1}v)l(v^{-1}s) dvdt 
\end{eqnarray*}
for any $f,g\in \cal L(B)$ and $k,l\in K(G)$. 
If we denote by $f\bullet k(r,s)=f(r)k(r^{-1}s)$, then 
$\delta_{\cal B}(f)(1\otimes k)(g\otimes l)(r,s) = 
(f\bullet k)(g\otimes l)(r,s)$. 
It is not hard to see that $\delta_{\cal B}$ is a full coaction
(note that $\delta_{\cal B}$ extends to a representation of 
${\cal L}_1({\cal B})$ and hence of $C^*({\cal B})$). 
\par\medskip

\begin{lemma} \label{1.6}
The set $N = \{ \sum_{i=1}^n f_i\bullet k_i:f_i\in {\cal L(B)}, 
k_i\in K(G)\}$ is dense in $C^*({\cal B}\times G)$. 
Consequently, $\delta_{\cal B}$ is a non-degenerate full coaction.
\end{lemma}
\noindent {\bf Proof}: 
It is sufficient to show that for any $f\in \cal L(B)$ and $k\in K(G)$, 
$f\otimes k$ can be approximated by elements in $N$ with respect 
to the ${\cal L}_1$-norm. 
Let $M$ and $K$ be the closed supports of $f$ and $k$ respectively. 
Since $k$ is uniformly continuous, for a given $\epsilon>0$, there exists 
a neighbourhood $V$ of $e$ such that $\|k(u)-k(v)\| < 
\epsilon/(\mu(M)\cdot \mu(MK)\cdot \sup_{t\in M}\|f(t)\|)$ if 
$u^{-1}v\in V$ (where $\mu$ is the left Haar measure of $G$).
Since $M$ is compact, there exist $r_1,..., r_n$ in $M$ such that 
$\cup_{i=1}^n r_iV$ covers $M$. 
Let $k_i(s) = k(r_i^{-1}s)$ and let $g_1,...,g_n$ be the partition of 
unity subordinate to $\{ r_iV\}_{i=1}^n$. 
Then
\begin{eqnarray*}
\| f\bullet k(r,s) - \sum_{i=1}^n (g_if\otimes k_i)(r,s)\| & \leq &
\|f(r)\| \cdot \mid \sum_{i=1}^n g_i(r)k(r^{-1}s) - 
\sum_{i=1}^n g_i(r)k(r_i^{-1}s)\mid \\
& \leq & \|f(r)\|\cdot (\ \sum_{i=1}^n\mid g_i(r)\mid\cdot 
\mid k(r^{-1}s) - k(r_i^{-1}s) \mid\ ). 
\end{eqnarray*}
As $g_i(r)\neq 0$ if and only if $r\in r_iV$, we have $\int_{MK}\int_M 
\|f\bullet k(r,s) - \sum_{i=1}^n (g_if\otimes k_i)(r,s)\| drds \leq 
\epsilon$. 
This proves the lemma.
\par\medskip

The following lemma is about general coactions of $C^*(G)$.
It implies, in particular that $\delta_\B$ is injective. 
Note that the trivial representation of $G$ on $\mathbb{C}$ 
induces a $*$-homomorphism $\cal E$ from $C^*(G)$ to $\mathbb{C}$ 
which is a coidentity in the sense that $({\cal E}\otimes \id)
\delta_G = \id = (\id\otimes {\cal E})\delta_G$ (where $\delta_G$ is 
the comultiplication on $C^*(G)$)
\par\medskip

\begin{lemma} \label{1.1} 
Let $\epsilon$ be a coaction on $A$ by $C^*(G)$. 
Suppose that $\epsilon_{\cal E} = (\id\otimes {\cal E})\epsilon$ and 
$A_{\cal E} = \epsilon_{\cal E} (A)$.
\par\noindent
(a) If $\epsilon$ is non-degenerate, then it is automatically 
injective. 
\par\noindent
(b) $A = A_{\cal E}\oplus \ker (\epsilon)$ (as Banach space). 
\end{lemma}
\noindent {\bf Proof}: 
(a) We first note that $\epsilon_{\cal E}$ is a 
$*$-homomorphism from $A$ to itself and so $A_{\cal E}$ is 
a $C^*$-subalgebra of $A$. 
It is clear that $\epsilon$ is injective on $A_{\cal E}$ and 
we want to show that $A=A_\epsilon$. 
For any $a\in A$ and any $s\in C^*(G)$ such that ${\cal E}(s)=1$, 
$a\otimes s$ can be approximated by elements of the form 
$\sum \epsilon(b_i)(1\otimes t_i)$ (as $\epsilon$ is non-degenerate). 
Therefore, $\sum \epsilon_{\cal E}(b_i){\cal E}(t_i) = 
(\id\otimes {\cal E})(\sum \epsilon(b_i)(1\otimes t_i))$ 
converges to $a$.
\par\noindent
(b) Note that for any $a\in A$, $\epsilon(a-\epsilon_{\cal E}(a)) = 0$. 
Since $\epsilon_{\cal E}$ is a projection on $A$, $A=A_{\cal E}\oplus 
\ker (\epsilon)$. 
\par\medskip

The above lemma actually holds for a general Hopf $C^*$-algebra with
a co-identity instead of $C^*(G)$ (for a brief review of Hopf 
$C^*$-algebras, please see e.g. [5] or [10]). 
\par\medskip

\begin{remark}
By [5, 7.15], if $\Gamma$ is a discrete amenable group, then any 
injective coaction of $C^*_r(\Gamma)$ is automatically non-degenerate. 
More generally, the arguments in [5, \S7] actually show that for 
any discrete group $G$, any injective coaction of $C^*(G)$ is 
non-degenerate. 
Hence, a coaction of $C^*(G)$ is injective if and only if it is 
non-degenerate (when $G$ is discrete). 
\end{remark}

We end this section with the following technical lemma.
\par\medskip

\begin{lemma} \label{1.3}
Let $A$ be a $C^*$-algebra and $E$ be a Hilbert $A$-module. 
Suppose that $({\cal H},\pi)$ is a faithful representation of $A$. 
Then 
\par\noindent 
(a) $\|x\| = \sup \{\|x\otimes_\pi \xi\|: \|\xi\|\leq 1\}$ for any $x\in E$; 
\par\noindent 
(b) the canonical map from ${\cal L}(E)$ to ${\cal L}(E\otimes_\pi {\cal H})$ 
(which sends $a$ to $a\otimes 1$) is injective. 
\end{lemma}

Part (a) follows from a direct computation and the
part (b) is a consequence of part (a). 
\par\medskip
\par\medskip
\par\medskip

\section{Reduced cross-sectional $C^*$-algebras}
\par\medskip
\par\medskip

In this section, we will define the reduced cross-sectional 
$C^*$-algebras for $C^*$-algebraic bundles and show that they 
carry canonical reduced coactions.
The intuitive idea is to consider the representation of 
${\cal L}_1({\cal B})$ as bounded operators on ${\cal L}_2({\cal B})$. 
However, since ${\cal L}_2({\cal B})$ is not a Hilbert $C^*$-module, it
seems unlikely that we can get a $C^*$-algebra out of this
representation. 
Instead, we will consider a slightly different version of 
``${\cal L}_2({\cal B})$'' which is a Hilbert $B_e$-module. 
The difficulty then is to show that the representation is well 
defined and bounded. 
This can be proved directly by a quite heavy analytical argument but we
will use Lemma \ref{2.8} to do the trick instead. 
We will also define the interesting notion of proper $C^*$-algebraic
bundles which will be needed in the next section. 
\par\medskip

\begin{lemma} \label{2.1} 
Consider the map $\langle \ ,\ \rangle _e$ from $\cal L(B)\times \cal L(B)$ 
to $B_e$ defined by
$$\langle f,g\rangle _e = \int_G f(t)^*g(t) dt$$ 
for all $f,g\in \cal L(B)$. 
Then $\langle \ ,\ \rangle _e$ is a $B_e$-valued inner product on 
$\cal L(B)$. 
\end{lemma}
\noindent {\bf Proof}: 
It is easily seen that $\langle \ ,\ \rangle_e$ is a well defined 
$B_e$-valued pre-inner product. 
Moreover, for all $f\in \cal L(B)$, $\langle f,f\rangle _e = 0$ 
if and only if $\int_G \varphi(f(t)^*f(t))dt = 0$ 
for all $\varphi\in (B_e)^*_+$ which implies that 
$f(t)^*f(t) = 0$ for all $t\in G$. 
\par\medskip

\begin{definition} \label{2.2} 
The completion of $\cal L(B)$ with respect to the $B_e$-valued inner 
product in Lemma \ref{2.1} is a Hilbert $B_e$-module and is denoted by 
\mbox{$(L^2_e({\cal B}), \|\cdot\|_e)$}. 
\end{definition}

It is clear that $\|\langle f,g\rangle _e\|\leq \|f\|_2\|g\|_2$ 
(by [7, II.5.4] and the H\" older's inequality). 
Hence there is a continuous map $J$ from ${\cal L}_2(\cal B)$ to 
$L^2_e({\cal B})$ with dense range. 
In fact, it is not hard to see that ${\cal L}_2(\cal B)$ is a right 
Banach $B_e$-module and $J$ is a module map. 
\par\medskip

Throughout this paper, $T$ is a non-degenerate *-representation 
of $\cal B$ on a Hilbert space ${\cal H}$ and $\phi$ is the 
restriction of $T$ on $B_e$. 
Moreover, $\mu_T$ is the representation of $C^*(\B)$ on $\cal H$ 
induced by $T$. 
By [7, VIII.9.4], $\phi$ is a non-degenerate representation of $B_e$. 
\par\medskip

\begin{lemma} \label{2.4}
There exists an isometry
$$V:L^2_e({\cal B})\otimes_\phi {\cal H}\rightarrow L^2(G; {\cal H}),$$
such that for all $f\in \cal L(B)$, $\xi\in\cal H$, and $t\in G$ one has
$$V(f\otimes \xi)(t) = T_{f(t)}\xi.$$
\end{lemma}
\noindent {\bf Proof}: 
It is easy to check that the map $V$ defined as above is inner 
product preserving and hence extends to the required map. 
\par\medskip

One technical difficulty in the study of reduced cross-sectional
$C^*$-algebras is that $V$ is not necessarily surjective. 
\par\medskip

\begin{example} \label{2.6} 
(a) If $\cal B$ is saturated, then $V$ is surjective. 
In fact, let $K= V(L^2_e({\cal B})\otimes_\phi {\cal H})$ and 
$\Theta$ be an element in the complement of $K$. 
For any $g\in \cal L(B)$ and $\eta\in {\cal H}$, 
$\int_G \langle T_{g(r)}\eta, \Theta(r)\rangle \ dr = 0$ which implies 
that $\int_G T_{g(r)}^*\Theta(r)\ dr = 0$.
Now for any $f\in \cal L(B)$, we have 
$$(\mu_T\otimes \lambda_G)\delta_{\cal B}(f)(\Theta)(t) = 
\int_G T_{f(s)}\Theta(s^{-1}t)\ ds = 
\int_G T_{f(tr^{-1})}\Theta(r)\Delta(r)^{-1}\ dr.$$ 
Moreover, for any $b\in B_{t^{-1}}$, let 
$g(r) = \Delta(r)^{-1}f(tr^{-1})^*b^*$. 
Then $g\in \cal L(B)$ and 
$$T_b (\mu_T\otimes \lambda_G)
\delta_{\cal B}(f)(\Theta)(t) = \int_G T_{g(r)}^*\Theta(r)\ dr = 0$$ 
for any $b\in B_{t^{-1}}$ (by the above equality). 
Since $\cal B$ is saturated and the restriction $\phi$ of $T$ is 
non-degenerate, $(\mu_T\otimes \lambda_G)\delta_{\cal B}(f)(\Theta)=0$
for any $f\in \cal L(B)$. 
Thus, $\Theta=0$ (because $(\mu_T\otimes \lambda_G)\circ\delta_{\cal B}$ 
is non-degenerate). 
\par\noindent 
(b) Let $\cal B$ be the trivial bundle over a discrete group $G$ 
(i.e. $B_e = \mathbb{C}$ and $B_t = (0)$ if $t\neq e$). 
Then $L^2_e({\cal B})\otimes_\phi {\cal H} \cong {\cal H}$ is a proper
subspace of $L^2(G; {\cal H})$.
\end{example}

For any $b\in B$, let $\hat T_b$ be the map from $L^2_e({\cal B})$ 
to itself defined by 
$\hat T_b(f) = b\cdot f$ for any $f\in \cal L(B)$ (where $b\cdot f(t) = 
bf(\pi(b)^{-1}t)$). 
The argument for $\hat T$ being continuous seems not easy. 
Instead, we will consider the corresponding representation of 
${\cal L}_1(\B)$ and show that it is well defined. 
\par\medskip

For any $f\in {\cal L}(\B)$, define a map $\lambda_\B (f)$ from ${\cal
L}(\B)$ to itself by $\lambda_\B (f)(g) = f\ast g$ ($g\in {\cal
L}(\B)$). 
We would like to show that this map is bounded and induces a bounded
representation of ${\cal L}_1(\B)$. 
In order to prove this, we will first consider a map 
$\tilde\lambda_\B(f)$ from ${\cal L}(\B)\otimes_{\rm alg} {\cal H}$
to itself given by 
$\tilde \lambda_\B(f)(g\otimes \xi) = f\ast g\otimes \xi$ ($g\in {\cal
L}(\B)$; $\xi\in \cal H$). 
In the following, we will not distinguish ${\cal L}(\B)\otimes_{\rm alg} 
{\cal H}$ and its image in $L_e^2(\B)\otimes_\phi{\cal H}$. 
\par\medskip

\begin{lemma} \label{2.8}
For any $f\in \cal L(B)$, $\tilde \lambda_{\cal B}(f)$ extends to a
bounded linear operator on $L_e^2(\B)\otimes_\phi{\cal H}$ such that
$\pil (f)\circ V = V\circ 
(\tilde \lambda_{\cal B}(f))$ 
(where $\pil$ is the composition:
$C^*(\B)\stackrel{\delta_\B}{\longrightarrow} C^*(\B)\otimes_{\rm max} C^*(G)
\stackrel{\mu_T\otimes \lambda_G}{\longrightarrow} \B({\cal H}\otimes
L^2(G))$). 
\end{lemma}
\noindent {\bf Proof:} 
For any $g\in \cal L(B)$, $\xi\in \cal H$ and $s\in G$, we have, 
\begin{eqnarray*}
\pil(f)V(g\otimes \xi)(s) 
&=& \int_G (T_{f(t)}\otimes \lambda_t) V(g\otimes \xi)(s) dt  
\ \ =\ \ \int_G T_{f(t)}T_{g(t^{-1}s)}\xi dt\\ 
&=& T_{f\ast g(s)}\xi
\ \ =\ \ V(f\ast g\otimes \xi)(s)
\ \ =\ \ V(\tilde \lambda_{\cal B}(f)(g\otimes\xi))(s).
\end{eqnarray*} 
Since $V$ is an isometry, $\tilde \lambda_{\cal B}(f)$ extends 
to a bounded linear operator on $L_e^2(\B)\otimes_\phi{\cal H}$ and 
satisfies the required equality. 
\par\medskip

Now by considering the representation $T$ for which $\phi$ is injective
and using Lemmas \ref{1.3}(a) and \ref{2.4}, $\lambda_\B(f)$ extends 
to a bounded linear map from $L_e^2(\B)$ to itself. 
It is not hard to show that $\langle f\ast g, h\rangle_e = 
\langle g,f^*\ast h\rangle_e$ (for any $g,h\in {\cal L}(\B)$). 
Hence $\lambda_\B(f)\in {\cal L}(L_e^2(\B))$.
Moreover, we have the following proposition. 
\par\medskip

\begin{proposition} \label{2.10}
The map $\lambda_{\cal B}$ from ${\cal L}_1(\cal B)$ to 
${\cal L}(L^2_e(\cal B))$ given by $\lambda_{\cal B} (f)(g) = f\ast g$ 
($f,g\in \cal L(B)$) is a well defined norm decreasing non-degenerate 
*-homomorphism such that $ \pil(f)\circ V = V\circ (\lambda_\B(f)
\otimes_\phi 1)$ ($f\in \cal L(B)$). 
\end{proposition}

\begin{definition} \label{2.11}
(a) $\lambda_\B$ is called the {\it reduced representation of $C^*(\B)$} 
and $C^*_r(\B) = \lambda_\B(C^*(\cal B))$ is called the 
{\it reduced cross-sectional $C^*$-algebra of $\cal B$}. 
\par\noindent 
(b) $\cal B$ is said to be {\it amenable} if $\lambda_{\cal B}$ is 
injective. 
\end{definition}

\begin{example} \label{2.12}
Suppose that $\cal B$ is the semi-direct product bundle corresponding to 
an action $\alpha$ of $G$ on a $C^*$-algebra $A$. 
Then $C^*({\cal B}) = A\times_\alpha G$ and $C^*_r({\cal B}) = 
A\times_{\alpha, r} G$. 
\end{example}

As in the case of full cross-sectional $C^*$-algebras, 
we can define non-degenerated reduced coactions on 
reduced cross-sectional $C^*$-algebras.
First of all, let us consider (as in the case of reduced 
group $C^*$-algebras) an operator $W$ from 
${\cal L(B}\times G)$ to itself defined by $W(F)(r,s) = 
F(r,r^{-1}s)$ ($F\in {\cal L(B}\times G)$). 
Note that for any $f\in\cal L(B)$ and $k\in K(G)$, $W(f\otimes k) = 
f\bullet k$ (where $f\bullet k$ is defined in the 
paragraph before Lemma \ref{1.6}) and that 
$L^2_e({\cal B}\times G) = L^2_e({\cal B})\otimes L^2(G)$ as 
Hilbert $B_e$-modules.
\par\medskip

\begin{lemma} \label{2.13} 
$W$ is a unitary in ${\cal L}(L^2_e({\cal B})\otimes L^2(G))$. 
\end{lemma}
\par
\noindent {\bf Proof}: 
For any $f,g\in \cal L(B)$ and $k,l\in K(G)$, we have the following
equality: 
\begin{eqnarray*}
\langle W(f\otimes k), W(g\otimes l)\rangle  
& = & \int_G\int_G (f\bullet k)(r,s)^*(g\bullet l)(r,s) dsdr\\ 
& = & \int_G\int_G f(r)^*\overline{k(r^{-1}s)} g(r)l(r^{-1}s) drds\\ 
& = & \int_G\int_G f(r)^*g(r)\overline{k(t)}l(t) dtdr 
\quad = \quad \langle f\otimes k, g\otimes l\rangle. 
\end{eqnarray*}
Hence $W$ is continuous and extends to an operator on 
$L^2_e({\cal B})\otimes L^2(G)$. 
Moreover, if we define $W^*$ by $W^*(f\otimes k)(r,s) = f(r)k(rs)$, 
then $W^*$ is the adjoint of $W$ and $WW^* = 1 = W^*W$.
\par\medskip

As in [12], we can define a *-homomorphism $\delta^r_{\cal B}$ 
from $C^*_r(\cal B)$ to ${\cal L}(L^2_e({\cal B})\otimes L^2(G))$ 
by $\delta^r_{\cal B}(x) = W(x\otimes 1)W^*$ 
($x\in C^*_r(\cal B)$). 
Moreover, for any $b\in B\subseteq M(C^*(\cal B))$ (see Lemma 
\ref{1.4}), $\delta^r_{\cal B}(\lambda_{\cal B}(b)) = 
\lambda_{\cal B}(b)\otimes \lambda_{\pi(b)}$ 
(where $\lambda_t$ is the canonical image of $t$ in $M(C^*_r(G))$).
\par\medskip

\begin{proposition} \label{2.14} 
The map $\delta^r_{\cal B}$ defined above is an injective non-degenerate 
coaction on $C^*_r(\cal B)$ by $C^*_r(G)$. 
\end{proposition}
\noindent {\bf Proof}: 
It is clear that $\delta^r_{\cal B}$ is an injective *-homomorphism. 
Moreover, 
$(\lambda_{\cal B}\otimes\lambda_G)\circ\delta_{\cal B} = 
\delta^r_{\cal B}\circ\lambda_{\cal B}$ 
which implies that $\delta^r_{\cal B}$ is a non-degenerate coaction
(see Lemma \ref{1.6}). 
\par\medskip

There is an alternative natural way to define ``reduced'' 
cross-sectional $C^*$-algebra (similar to the corresponding 
situation of full and reduced crossed products): 
$C^*_R({\cal B}) := C^*({\cal B})/\ker (\epsilon_{\cal B})$ 
(where $\epsilon_{\cal B}$ is the composition: $C^*(\B)\stackrel
{\delta_{\cal B}}{\longrightarrow}C^*(\B)\otimes_{\rm max}C^*(G)
\stackrel{\id\otimes \lambda_G}{\longrightarrow}C^*(\B)\otimes C^*_r(G)$).
\par\medskip

\begin{remark} \label{2.b}
(a) It is clear that $\pil = (\mu_T\otimes \lambda_G)\circ\delta_\B$
(see Lemma \ref{2.8}) induces a representation of $C^*_R(\B)$ on 
$L^2(G; {\cal H})$. 
If $\mu_T$ is faithful, then this induced representation is also
faithful and $C^*_R(\B)$ can be identified with the image of $\pil$. 
\par\noindent
(b) If $\mu_T$ is faithful, then so is $\phi$ and 
$\lambda_\B\otimes_\phi 1$
is a faithful representation of $C^*_r(\B)$ (by Lemma \ref{1.3}(b)). 
Therefore, part (a) and Lemma \ref{2.10} implies that $C^*_r(\B)$ is a
quotient of $C^*_R(\B)$. 
\end{remark}

In [12, 3.2(1)], it was proved that these two reduced cross-sectional
$C^*$-algebras coincide in the case of semi-direct product bundles. 
The corresponding result in the case of $C^*$-algebraic bundles 
over discrete groups was proved implicitly in [6, 4.3]. 
In the following we shall see that it is true in general. 
\par\medskip

The idea is to define a map $\varphi$ from $C_r^*(\cal B)$ to
${\cal L}(L^2(G; {\cal H}))$ such that $\varphi\circ\lambda_{\cal B} =
\pil$ (see Remark \ref{2.b}(a)). 
As noted above, the difficulty is that $V$ may not be surjective and,
by Lemma 2.10, $\lambda_{\cal B}\otimes_\phi 1$ may only be a proper 
subrepresentation of $\pil$ (see Example \ref{2.6}(b)). 
However, we may ``move it around'' filling out the whole 
representation space for $\pil$ using the right regular
representation $\rho$ of $G$ (on $L^2(G)$):  
$\rho_r(g)(s)=\Delta(r)^{1/2}g(sr)$ ($g\in L^2(G)$; $r,s\in G$) where
$\Delta$ is the modular function for $G$.
\par\medskip

\begin{lemma} \label{2.a}
For each $r\in G$,
\par\noindent
  (a) The unitary operator $\rho_r\otimes 1$ on $L_2(G)\otimes {\cal H}= 
\LGH$ lies in the commutant of $\pil (C^*(\cal B))$.
\par\noindent
  (b) Consider the isometry
  $$
V^r: L_2(\B)\otimes_\phi{\cal H} \to \LGH,
  $$
given by $V^r = (\rho_r\otimes 1) V$.  
Then for all $a\in C^*(\B)$ one
has $V^r (\lambda_{\cal B}(a)\otimes 1) = \pil(a) V^r$.
\par\noindent
  (c) Let $K_r$ be the range of\/ $V^r$.  
Then $K_r$ is invariant under $\pil$ and the restriction of $\pil$ 
to $K_r$ is equivalent to $\lambda_{\cal B}\otimes 1$.
\end{lemma}
\noindent {\bf Proof:} 
It is clear that $\rho_r\otimes 1$ commutes with $\pil(b_t) =
\lambda_t\otimes T_{b_t}$ for any $b_t\in B_t$ (see Lemma \ref{1.4}).
It then follows that $\rho_r\otimes 1$ also commutes with the 
range of the integrated form of $\pil$, whence (i).  
The second point follows immediately from (i) and
Proposition \ref{2.10}. Finally, (iii) follows from (ii).
\par\medskip

  Our next result is intended to show that the $K_r$'s do indeed fill
out the whole of $\LGH$.
\par\medskip

\begin{proposition} \label{2.c}
  The linear span of\/ $\bigcup_{r\in G} K_r$ is dense in $\LGH$.
\end{proposition}
\noindent {\bf Proof:} 
Let
  $$
\Gamma = {\rm span}\{V^r( f\otimes \eta): r\in G,\  f\in {\cal L}(\B),
\ \eta\in{\cal H}\}.
  $$
Since for any $t\in G$,
  $$
V^r( f\otimes \eta)(t) =
(\rho_r\otimes 1)V( f\otimes \eta)(t) =
\Delta(r)^{1/2}V(f\otimes \eta)(tr) =
\Delta(r)^{1/2}T_{f(tr)}\eta,
  $$
and since we are taking $f$ in ${\cal L}(\B)$ above, it is easy to see
that $\Gamma$ is a subset of $C_c(G,{\cal H})$.  
Our strategy will be to use [7, II.15.10] 
(on the Banach bundle ${\cal H}\times G$ over $G$) 
for which we must prove that:
\par
\noindent (I) If $f$ is a continuous complex function on $G$ and
$\zeta\in\Gamma$, then the pointwise product $f\zeta$ is in $\Gamma$;
\par
\noindent (II) For each $t\in G$, the set $\{\zeta(t):\zeta\in \Gamma\}$
is dense in ${\cal H}$.
\par\noindent
The proof of (I) is elementary in view of the fact that ${\cal L}(\B)$ is
closed under pointwise multiplication by continuous scalar-valued
functions [7, II.13.14].
  In order to prove (II), let $\xi\in{\cal H}$ have the form
$\xi=T_b\eta=\phi(b)\eta$, where $b\in B_e$ and $\eta\in{\cal H}$.  
  By [7, II.13.19], let $f\in {\cal L}(\B)$ be such that $f(e)=b$.
  It follows that
$\zeta_r:=V^r(f\otimes \eta)$ is in $\Gamma$ for all $r$.
  Also note that,
setting $r=t^{-1}$, we have
  $$
\zeta_{t^{-1}}(t) =
\Delta(t)^{-1/2}T_{ f(e)}\eta =
\Delta(t)^{-1/2}T_{b}\eta =
\Delta(t)^{-1/2}\xi.
  $$
  This shows that $\xi\in\{\zeta(t):\zeta\in \Gamma\}$.  
  Since the set of such $\xi$'s is dense in ${\cal H}$ (because 
$\phi$ is non-degenerate by assumption), we have that (II) is proven.
  As already indicated, it now follows from [7, II.15.10] that
$\Gamma$ is dense in $\LGH$.  Since $\Gamma$ is contained in the
linear span of $\bigcup_{r\in G} K_r$, the conclusion follows.
\par\medskip

  We can now obtain the desired result. 
\par\medskip

\begin{theorem} \label{2.d}
For all $a\in C^*(\B)$ one has that $\|\pil(a)\| \leq
\|\lambda_{\cal B}(a)\|$.
Consequently, $C^*_R(\B)=C^*_r(\B)$. 
\end{theorem}
\noindent {\bf Proof:}
  We first claim that for all $a\in C^*(\B)$ one has
that
  $$
  \lambda_{\cal B}(a) = 0 \quad =\!\!\Rightarrow \quad \pil(a) = 0.
  $$
  Suppose that $\lambda_{\cal B}(a) = 0$.
  Then for each $r\in G$ we have by Lemma \ref{2.a}(b) that
  $$
  \pil(a) V^r = V^r (\lambda_{\cal B}(a)\otimes 1) =0.
  $$
  Therefore
  $\pil(a)=0$ in the range $K_r$ of $V^r$.
  By Proposition \ref{2.c} it follows that $\pil(a)=0$, thus proving our
claim.
  Now define a map
  $$
  \varphi : C^*_r(\B) \longrightarrow \B(\LGH)
  $$
by $\varphi (\lambda_{\cal B}(a)) := \pil(a)$, for all $a$ in $C^*(\B)$.  
  By the claim above we have that $\varphi$ is well defined.  
  Also, it is easy to see that $\varphi$ is a *-homomorphism.  
  It follows that $\varphi$ is contractive and hence that for all 
$a$ in $C^*(\B)$
  $$
\|\pil(a)\| =
\|\varphi (\lambda_{\cal B}(a))\| \leq
\|\lambda_{\cal B}(a)\|.
  $$
  For the final statement, we note that if $\mu_T$ is faithful, then 
the map $\varphi$ defined above is the inverse of the quotient map 
from $C^*_R(\B)$ to $C^*_r(\B)$ given in Remark \ref{2.b}(ii). 
\par\medskip

  The following generalises [11, 7.7.5] to the context of $C^*$-
algebraic bundles:
\par\medskip

\begin{corollary} \label{2.e}
  Let $T:\B\rightarrow {\cal L}({\cal H})$ be a non-degenerate 
$*$-representation of the $C^*$-algebraic bundle $\B$
and let $\pil$ be the representation of $\B$ on $\LGH$ given by
$\pil(b_t) = \lambda_t\otimes T_{b_t}$, for $t\in G$, and $b_t\in B_t$.
  Then $\pil$ is a well defined representation and induces a 
representation of $C^*(\B) ($again denoted by $\pil$).
  In this case, $\pil$ factors through $C^*_r(\B)$.
  Moreover, if $\phi = T\!\mid_{B_e}$ is faithful, the 
representation of $C^*_r(\B)$ arising from this factorisation 
is also faithful.
\end{corollary}
\noindent {\bf Proof:}
By Remark \ref{2.b}(a), $\pil$ factors through a representation of 
$C^*_R(\B)=C^*_r(\B)$ (Theorem \ref{2.d}). 
Now if $\phi$ is faithful, then by Theorem \ref{2.d}, Lemmas \ref{2.10}
and \ref{1.3}(a), $\|\pil(a)\| = \|\lambda_{\cal B}(a)\|$. 
This proves the second statement. 
\par\medskip
\par\medskip
\par\medskip

\section{The approximation property of $C^*$-algebraic bundles}
\par\medskip
\par\medskip

From now on, we assume that $\mu_T$ (see the paragraph 
before Lemma \ref{2.4}) is faithful.
Moreover, we will not distinguish ${\cal L}(\B)$ and its image in $C^*(\B)$. 
\par\medskip 

The materials in this section is similar to the discrete case in [6].
Let ${\cal B}_e$ be the $C^*$-algebraic bundle $B_e\times G$ over $G$.
We will first define a map from $L^2_e({\cal B}_e)\times 
C^*_r({\cal B})\times L^2_e({\cal B}_e)$ to $C^*(\cal B)$. 
For any $\alpha\in {\cal L(B}_e)$, let $V_\alpha$ be a map 
from ${\cal H}$ to $L^2(G; {\cal H})$ given by
$$V_\alpha(\xi)(s) = \phi(\alpha(s))\xi$$ 
($\xi\in {\cal H}$; $s\in G$). 
It is clear that $V_\alpha$ is continuous and $\|V_\alpha\|\leq \|\alpha\|$. 
Moreover, we have $V_\alpha^*(\Theta) = 
\int_G \phi(\alpha(r)^*)\Theta(r)\ dr$ 
($\alpha\in {\cal L(B}_e)$; $\Theta\in L^2(G)\otimes {\cal H} = 
L^2(G;{\cal H})$) and $\|V_\alpha^*\|\leq \|\alpha\|$. 
Thus, for any $\alpha, \beta\in L^2_e({\cal B}_e)$, we obtain a continuous 
linear map $\Psi_{\alpha,\beta}$ from ${\cal L}(L^2(G)\otimes {\cal H})$ 
to ${\cal L}({\cal H})$ defined by
$$\Psi_{\alpha,\beta}(x) = V_\alpha^*xV_\beta$$ 
with $\|\Psi_{\alpha,\beta}\|\leq \|\alpha\|\|\beta\|$. 
Recall from Remark \ref{2.b}(a) and Theorem \ref{2.d} that $C^*_r(\cal B)$ 
is isomorphic to the image of $C^*(\cal B)$ in ${\cal L}(L^2(G)\otimes 
{\cal H})$ under $\pil = (\mu_T\otimes \lambda_G)\circ \delta_{\cal B}$. 
\par\medskip

\begin{lemma} \label{3.1}
Let $\alpha, \beta\in {\cal L(B}_e)$ and $f\in \cal L(B)$. 
Then $\Psi_{\alpha,\beta}(\pil(f)) = 
\alpha\cdot f\cdot \beta$ where $\alpha\cdot f\cdot \beta\in {\cal L(B)}$
is defined by $\alpha\cdot f\cdot \beta(s) = 
\int_G \alpha(t)^*f(s)\beta(s^{-1}t)\ dt$. 
\end{lemma}
\noindent {\bf Proof}: 
For any $\xi \in {\cal H}$, we have
\begin{eqnarray*}
\Psi_{\alpha,\beta}(\pil(f))\xi & = & 
\int_G \phi(\alpha(t)^*)(\pil(f)V_\beta\xi)(t)\ dt\\ 
& = & \int_G\int_G \phi(\alpha^*(t))T_{f(s)}\phi(\beta(s^{-1}t))
\xi\ ds dt\\ 
& = & \int_G T_{(\alpha\cdot f\cdot \beta)(s)}\xi\ ds. 
\end{eqnarray*}
\par\medskip

Hence we have a map from $L^2_e({\cal B}_e)\times 
C^*_r({\cal B})\times L^2_e({\cal B}_e)$ to $C^*(\cal B)$ such 
that $\|\alpha\cdot x\cdot \beta\|\leq \|\alpha\|\|x\|\|\beta\|$. 
Next, we will show that this map sends 
$L^2_e({\cal B}_e)\times \pil({\cal L(B}))\times 
L^2_e({\cal B}_e)$ to $\cal L(B)$. 
\par\medskip

\begin{lemma} \label{3.2}
For any $\alpha, \beta\in L^2_e({\cal B}_e)$ and $f\in \cal L(B)$, 
$\alpha\cdot f\cdot \beta\in \cal L(B)$. 
\end{lemma}
\noindent {\bf Proof}: 
If $\alpha', \beta' \in {\cal L(B}_e)$
\begin{eqnarray*}
\|(\alpha'\cdot f\cdot \beta')(s)\| & = &
\sup \{\mid \langle \eta,\int_G T_{\alpha'(t)^*f(s)\beta'(s^{-1}t)}
\xi\ dt\rangle \mid: \|\eta\|\leq 1; \|\xi\|\leq 1\}\\ 
& \leq & \sup \{ \|f(s)\| (\int_G \|\phi(\alpha'(t))\eta\|^2\ dt)^{1/2} 
(\int_G \|\phi(\beta'(t))\xi\|^2\ dt)^{1/2}: \|\eta\|\leq 1; \|\xi\|\leq
1\}\\ 
& = & \|f(s)\|\|\alpha'\|\|\beta'\|. 
\end{eqnarray*}
Let $\alpha_n$ and $\beta_n$ be two sequences of elements in ${\cal L(B}_e)$ 
that converge to $\alpha$ and $\beta$ respectively. 
Then $(\alpha_n\cdot f\cdot \beta_n)(s)$ converges to an element 
$g(s)\in B_s$. 
Moreover, since $f$ is of compact support and the convergence is uniform, 
$g\in \cal L(B)$ and ${\rm supp}(g)\subseteq {\rm supp}(f)$. 
In fact, this convergence actually takes place in 
${\cal L}_1({\cal B})$ and hence in $C^*(\cal B)$. 
Therefore $\Psi_{\alpha,\beta}(\pil(f))=\mu_T(\alpha\cdot f\cdot \beta)$. 
\par\medskip

\begin{remark} \label{3.3}
The proof of the above lemma also shows that $\Psi_{\alpha,\beta}$ 
sends the image of $B_t$ in $M(C^*_r(\cal B))$ to the image 
of $B_t$ in $M(C^*(\cal B))$. 
Hence $\Psi_{\alpha,\beta}$ induces a map $\Phi_{\alpha,\beta}$ from $B$
to $B$ which preserves fibers. 
\end{remark}

\begin{definition} \label{3.4} 
Let $\{\Phi_i\}$ be a net of maps from $B$ to itself such that 
they preserve fibers and are linear on each fiber.
\par\noindent
(a) $\{\Phi_i\}$ is said to be {\it converging to $1$ on compact slices 
of $B$} if for any $f\in {\cal L(B})$ and any $\epsilon > 0$, there
exists $i_0$ such that for any $i\geq i_0$, $\|\Phi_i(b)-b\|<\epsilon$
for any $b\in f(G)$ ($f(G)$ is called a {\it compact slice} of B). 
\par\noindent 
(b) Then $\{\Phi_i\}$ is said to be {\it converging to 
$1$ uniformly on compact-bounded subsets of $B$} if for any compact 
subset $K$ of $G$ and any $\epsilon > 0$, there exists $i_0$ such that 
for any $i\geq i_0$, $\|\Phi_i(b)-b\|<\epsilon$ if 
$\pi(b)\in K$ and $\|b\|\leq 1$. 
\end{definition}

\begin{lemma} \label{3.5}
Let $\{\Phi_i\}$ be a net as in Definition \ref{3.4}.
Then each of the following conditions is stronger than the next one. 
\begin{enumerate}
\item[i.] $\{\Phi_i\}$ converges to $1$ uniformly on 
compact-bounded subsets of $B$. 
\item[ii.] $\{\Phi_i\}$ converges to $1$ uniformly on 
compact slices of $B$.
\item[iii.] For any $f\in \cal L(B)$, the net $\Phi_i\circ f$ converges 
to $f$ in ${\cal L}_1(\cal B)$. 
\end{enumerate}
\end{lemma}
\par\noindent {\bf Proof:}
Since every element in ${\cal L(B})$ has compact support and is bounded,
it is clear that (i) implies (ii). 
On the other hand, (ii) implies (iii) is obvious.
\par\medskip

Following the idea of [6], we define the approximation property 
of $\cal B$. 
\par\medskip

\begin{definition} \label{3.6}
(a) Let $\cal B$ be a $C^*$-algebraic bundle.
For $M>0$, $\cal B$ is said to have the {\it $M$-approximation property} 
(respectively, {\it strong $M$-approximation property}) if 
there exist nets $(\alpha_i)$ and $(\beta_i)$ in ${\cal L(B}_e)$ such that 
\begin{enumerate} 
\item[i.] $\sup_i \|\alpha_i\|\|\beta_i\| \leq M$;
\item[ii.] the map $\Phi_i = \Phi_{\alpha_i,\beta_i}$ (see Remark \ref{3.3}) 
converges to $1$ uniformly on compact slices of $B$ (respectively, 
uniformly on compact-bounded subsets of $B$). 
\end{enumerate} 
$\cal B$ is said to have the (respectively, {\it strong}) 
{\it approximation property} if it has 
the (respectively, strong) $M$-approximation property 
for some $M > 0$.
\par\noindent 
(b) We will use the terms ({\it strong}) {\it positive 
$M$-approximation property} and ({\it strong}) {\it positive 
approximation property} if we can choose $\beta_i = \alpha_i$ 
in part (a). 
\end{definition}

Because of Remark 3.7(b) below as well as [11, 7.3.8], we believe that
the above is the weakest condition one can think of to ensure the
amenability of the $C^*$-algebraic bundle. 
\par\medskip

\begin{remark} \label{3.7}
(a) Since any compact subset of a discrete group is finite and any 
$C^*$-algebraic bundle has enough cross-sections, 
the approximation property defined in [6] is the
same as positive 1-approximation property defined above. 
\par\noindent
(b) It is easy to see that the amenability of $G$ implies 
the positive 1-approximation property of $\cal B$ (note that the 
positive 1-approximation property is similar to the condition in 
[11, 7.3.8]). 
In fact, let $\xi_i$ be the net given by [11, 7.3.8] and let 
$\eta_i(t) = \overline{\xi_i(t)}$. 
If $u_j$ is an approximate unit of $B_e$ (which is also a strong 
approximate unit of $\cal B$ by [7, VIII.16.3]), 
then the net $\alpha_{i,j} = \beta_{i,j} = \eta_i u_j$
will satisfy the required property.
\par\noindent
(c) We can also formulate the approximation property as follows: 
there exists $M>0$ such that for any compact slice $S$ of $B$ and any 
$\epsilon > 0$, there exist $\alpha,\beta\in L^2_e({\cal B}_e)$ with 
$$\|\alpha\|\|\beta\|\leq M\qquad {\rm and} \qquad 
\|\alpha\cdot b\cdot \beta - b\| < \epsilon$$ 
if $b\in S$. 
In fact, we can replace $L^2_e({\cal B}_e)$ by ${\cal L(B}_e)$ and 
consider the directed set $D= \{ (K,\epsilon): K$ is a compact subset 
of $G$ and $\epsilon > 0 \}$. 
For any $d=(K,\epsilon)\in D$, we take $\alpha_d$ and $\beta_d$ 
that satisfying the above condition. 
These are the required nets. 
\end{remark}

We can now prove the main results of this section. 
\par\medskip

\begin{proposition} \label{3.8}
If $\cal B$ has the approximation property, then the coaction 
$\epsilon_{\cal B} = (\id\otimes \lambda_G)\circ\delta_{\cal B}$ is injective. 
\end{proposition}
\noindent {\bf Proof}: 
Let $\Phi_i = \Phi_{\alpha_i,\beta_i}$ be the map from $B$ to itself 
as given by Definition \ref{3.6}(a)(ii) and $\Psi_i = \Psi_{\alpha_i,\beta_i}$. 
Let $J_i=\Psi_i\circ\pil$. 
By Lemma \ref{3.2}, for any $f\in\cal L(B)$, $J_i(f)\in {\cal L}(\B)$ 
(note that we regard ${\cal L}(\B)\subseteq C^*(\B)$) and $J_i(f)(s) = 
\Phi_i(f(s))$ ($s\in G$).
Since $\Phi_i\circ f$ converges to $f$ in ${\cal L}_1(\cal B)$ 
(by Lemma \ref{3.5}), $J_i(f)$ converges to $f$ in $C^*(\cal B)$. 
Now because $\|J_i\|\leq \|\Psi_i\|\leq \sup_i\|\alpha_i\|
\|\beta_i\|\leq M$, we know that $J_i(x)$ converges to $x$ for all 
$x\in C^*(\cal B)$ and $\epsilon_{\cal B}$ is injective. 
\par\medskip

Note that if $G$ is amenable, we can also obtain directly from 
the Lemma \ref{1.1}(a) that $\epsilon_{\cal B}$ is injective. 
\par\medskip

\begin{theorem} \label{3.9}
Let $\cal B$ be a $C^*$-algebraic bundle having the 
approximation property (in particular, if $G$ is amenable). 
Then $\cal B$ is amenable. 
\end{theorem}
\noindent {\bf Proof:}
Proposition \ref{3.8} implies that $C_R^*(\B) = C^*(\B)$
(see the paragraph before Remark \ref{2.b}). 
Now the amenability of $\B$ clearly follows from Theorem \ref{2.d}.
\par\medskip
\par\medskip
\par\medskip

\section{Two special cases}
\par\medskip
\par\medskip

\noindent {\em I. Semi-direct product bundles and nuclearity of 
crossed products.}
\par\medskip

Let $A$ be a $C^*$-algebra with action $\alpha$ by a locally compact 
group $G$. 
Let ${\cal B}$ be the semi-direct product bundle of $\alpha$. 
\par\medskip

\begin{remark} \label{4.1}
$\cal B$ has the (respectively, strong) $M$-approximation property 
if there exist nets $\{\gamma_i\}_{i\in I}$ and 
$\{\theta_i\}_{i\in I}$ in $K(G;A)$ such that 
$$\|\int_G \gamma_i(r)^*\gamma_i(r)\ 
dr\|\cdot\|\int_G \theta_i(r)^*\theta_i(r)\ dr\| \leq M^2$$
and for any $f\in K(G;A)$ (respectively, for any compact subset 
$K$ of $G$), 
$\int_G \gamma_i(r)^*a\alpha_t (\theta_i(t^{-1}r))\ dr$ converges to
$a\in A$ uniformly for $(t,a)$ in the graph of $f$ (respectively,
uniformly for $t\in K$ and $\|a\|\leq 1$). 
\end{remark}

\begin{definition} \label{4.2} 
An action $\alpha$ is said to have the (respectively, strong) 
{\it ($M$-)approximation property} 
(respectively, $\alpha$ is said to be {\it weakly amenable}) if 
the $C^*$-algebraic bundle $\cal B$ associated with $\alpha$ has the 
(respectively, strong) ($M$-)approximation property 
(respectively, $\B$ is amenable). 
\end{definition}

Let $G$ and $H$ be two locally compact groups. 
Let $A$ and $B$ be $C^*$-algebras with actions $\alpha$ and $\beta$ by 
$G$ and $H$ respectively. 
Suppose that $\tau = \alpha\otimes \beta$ is the product action on 
$A\otimes B$ by $G\times H$. 
\par\medskip

\begin{lemma} \label{4.3} 
With the notation as above, if $A$ is nuclear and both $\alpha$ and 
$\beta$ have the approximation property, then $(A\otimes B)
\times_\tau(G\times H)=(A\otimes B)\times_{\tau,r}(G\times H)$. 
\end{lemma}
\noindent {\bf Proof}: 
Let $\cal B$, $\cal D$ and $\cal F$ be the semi-direct product bundles of 
$\alpha$, $\beta$ and $\tau$ respectively. 
Then $C^*_r({\cal F})=C^*_r({\cal B})\otimes C^*_r({\cal D})$ 
(by Example \ref{2.12}). 
Moreover, since $A$ is nuclear, $C^*({\cal F}) = 
C^*({\cal B})\otimes_{\max} C^*({\cal D})$ (by Example \ref{2.12} and 
[10, 3.2]). 
It is not hard to see that the coaction, $\delta_{\cal F}$, on 
$C^*(\cal F)$ is the tensor product of the coactions on 
$C^*(\cal B)$ and $C^*(\cal D)$. 
Suppose that $C^*(\cal F)$ is a $C^*$-subalgebra of $\cal L(H)$.
Consider as in Section 2, the composition:
$$\mu_{\cal F}: C^*({\cal F})\stackrel{\delta_{\cal F}}
{\longrightarrow} C^*({\cal F})\otimes_{\rm max} C^*(G\times H)
\stackrel{\id\otimes \lambda_{G\times H}}{\longrightarrow}
{\cal L} ({\cal H}\otimes L^2(G\times H))$$ 
and identify its image with $C^*_R({\cal F})=C^*_r({\cal F})$ 
(see Remark \ref{2.b}(a)).
We also consider similarly the maps $\mu_{\cal B}$ and $\mu_{\cal D}$ 
from $C^*(\cal B)$ and $C^*(\cal D)$ to ${\cal L}({\cal H}\otimes
L^2(G))$ and ${\cal L}({\cal H}\otimes L^2(H))$ respectively. 
Now for any $f\in \cal L(B)$ and $g\in \cal L(D)$, 
we have $\mu_{\cal F}(f\otimes g) = (\mu_{\cal B}(f))_{12} 
(\mu_{\cal D}(g))_{13}\in {\cal L}({\cal H}\otimes L^2(G)\otimes L^2(H))$. 
As in Section 3, we define, for any $k\in {\cal L(B}_e)$ and 
$l\in {\cal L(D}_e)$, an operator $V_{k\otimes l}$ from $\cal H$ to 
${\cal H}\otimes L^2(G\times H)$ by 
$V_{k\otimes l}\zeta (r,s) = k(r)(l(s)\zeta)$
($r\in G$; $s\in H$; $\zeta\in {\cal H}$). 
It is not hard to see that $V_{k\otimes l}(\zeta) = 
(V_k\otimes 1)V_l(\zeta)$ and
$$V_{k\otimes l}^*\mu_{\cal F}
(f\otimes g)V_{k'\otimes l'} = (V_k^*\mu_{\cal B}(f)V_{k'})
(V_l^*\mu_{\cal D}(g)V_{l'})\in \cal L(H)$$ 
(note that $B_r$ 
commutes with $D_s$ in ${\cal L}(\cal H)$). 
Now let $k_i, k'_i\in \cal L(B)$ and $l_j, l'_j\in \cal L(D)$ be 
the nets that give the corresponding approximation property on 
$\cal B$ and $\cal D$ respectively. 
Then $V_{k_i}^*\mu_{\cal B}(f)V_{k'_i}$ converges to $f$ in 
$C^*(\cal B)$ and $V_{l_j}^*\mu_{\cal D}(g)V_{l'_j}$ 
converges to $g$ in $C^*(\cal D)$. 
Hence $J_{i,j}(z) = 
V_{k_i\otimes l_j}^*\mu_{\cal F}(z)V_{k'_i\otimes l'_j}$ 
converges to $z$ in $C^*(\cal F)$ 
for all $z\in {\cal L(B)}\otimes_{alg} {\cal L(D)}$. 
Since $\|J_{i,j}\|$ is uniformly bounded, $\mu_{\cal F}$ is 
injective and $C^*({\cal F}) = C^*_r({\cal F})$. 
\par\medskip

An interesting consequence of this proposition is the nuclearity 
of the crossed products of group actions with the approximation 
property (which is a generalisation of the case of actions 
of amenable groups). 
Note that in the case of discrete group, this was also proved by 
Abadie in [1].
\par\medskip

\begin{theorem} \label{4.4} 
Let $A$ be a $C^*$-algebra and $\alpha$ be an action on $A$ by a 
locally compact group $G$. 
If $A$ is nuclear and $\alpha$ has the approximation property, 
then $A\times_{\alpha} G = A\times_{\alpha , r} G$ 
is also nuclear. 
\end{theorem}
\noindent {\bf Proof}: 
By Lemma \ref{4.3} (or Theorem \ref{3.9}), 
$A\times_{\alpha} G = A\times_{\alpha , r} G$.
For any $C^*$-algebra $B$, let $\beta$ be the trivial action
on $B$ by the trivial group $\{e\}$.
Then $(A\times_\alpha G)\otimes_{\max} B =
(A\otimes_{\max}B)\times_{\alpha\otimes\beta} G = 
(A\otimes B)\times_{\alpha\otimes\beta, r} G = 
(A\times_{\alpha , r} G)\otimes B$ (by Lemma \ref{4.3} again). 
\par\medskip

One application of Theorem \ref{4.4} is to relate the amenability of 
Anantharaman-Delaroche (see [4, 4.1]) to the approximation property 
in the case when $A$ is nuclear and $G$ is discrete. 
The following corollary clearly follows from this theorem and [4, 4.5]. 
\par\medskip

\begin{corollary} \label{4.5}
Let $A$ be a nuclear $C^*$-algebra with an action $\alpha$ by 
a discrete group $G$. 
If $\alpha$ has the approximation property, 
then $\alpha$ is amenable in the sense of Anantharaman-Delaroche. 
\end{corollary}

We don't know if the two properties coincide in general. 
However, we can have a more direct and transparent comparison of 
them in the case of commutative $C^*$-algebras and show that they 
are the same.
Furthermore, they also coincide in the case of finite 
dimensional $C^*$-algebras. 
\par\medskip

\begin{corollary} \label{4.a}
Let $A$ be $C^*$-algebra with action $\alpha$ by a discrete group $G$. 
\par\noindent
(a) If $A$ is commutative, the followings are equivalent: 
\begin{enumerate}
\item[i.] $\alpha$ is amenable in the sense of Anantharaman-Delaroche;
\item[ii.] $\alpha$ has the positive 1-approximation property;
\item[iii.] $\alpha$ has the approximation property. 
\end{enumerate}
\par\noindent
(b) If $A$ is unital and commutative or if $A$ is finite 
dimensional, then (i)-(iii) are also equivalent to the 
following conditions: 
\begin{enumerate}
\item[iv.] $\alpha$ has the strong positive 1-approximation property;
\item[v.] $\alpha$ has the strong approximation property. 
\end{enumerate}
\end{corollary}
\noindent {\bf Proof:}
(a) By [4, 4.9(h')], $\alpha$ is amenable in the sense of 
Anantharaman-Delaroche if and only if 
there exists a net $\{\gamma_i\}$ in $K(G;A)$ such that 
$\|\sum_{r\in G} \gamma_i(r)^*\gamma_i(r)\|\leq 1$ 
and $\sum_{r\in G} \gamma_i(r)^*\alpha_t(\gamma_i(t^{-1}r))$ converges 
to $1$ strictly for any $t\in G$. 
It is exactly the original definition of the approximation property
given in [1].
Hence conditions (i) is equivalent to conditions (ii) (see Remark 
\ref{3.7}(a)). 
Now part (a) follows from Corollary \ref{4.5}. 
\par\noindent
(b) Suppose that $A$ is both unital and commutative. 
Let $\alpha$ satisfy condition (i) and $\{\gamma_i\}$ be the net as
given in the proof of part (a) above. 
As $A$ is unital, the strict convergence and the norm convergence are
equivalent. 
Moreover, as $G$ is discrete, any compact subset $K$ of $G$ is finite. 
These, together with the commutativity of $A$, imply that 
$\sum_{r\in G} \gamma_i(r)^*\alpha_t(\gamma_i(t^{-1}r))$ converges to 1
strictly for any $t\in G$ if and only if 
$\sum_{r\in G} \gamma_i(r)^*a\alpha_t(\gamma_i(t^{-1}r))$ converges to 
$a\in A$ uniformly for $t\in K$ and $\|a\|\leq 1$. 
Thus, by Remark \ref{4.1}, we have the equivalence of (i) and (iv) 
in the case of commutative unital $C^*$-algebras and the equivalence 
of (i)-(v) follows from Lemma \ref{3.5} and Corollary \ref{4.5}.
Now suppose that $A$ is a finite dimensional $C^*$-algebra (but not 
necessary commutative). 
By [4, 4.1] and [4, 3.3(b)], $\alpha$ satisfies condition (i) if and 
only if there exists a net $\{\gamma_i\}$ in $K(G; Z(A))$ 
(where $Z(A)$ is the centre of $A=A^{**}$) such that 
$\|\sum_{r\in G} \gamma_i(r)^*\gamma_i(r)\|\leq 1$ 
and for any $t\in G$, $\sum_{r\in G} \gamma_i(r)^*\alpha_t
(\gamma_i(t^{-1}r))$ converges to $1$ weakly (and hence converges 
to $1$ in norm as $A$ is finite dimensional). 
Let $K$ be any compact (and hence finite) subset of $G$. 
Since $\alpha_t(\gamma_i(t^{-1}r))\in Z(A)$, 
$\sum_{r\in G} \gamma_i(r)^*a \alpha_t(\gamma_i(t^{-1}r))$ 
converges to $a\in A$ uniformly for $t\in K$ and $\|a\|\leq 1$. 
This shows that $\alpha$ satisfies condition (iv) 
(see Remark \ref{4.1}).
The equivalence follows again from Lemma \ref{3.5} and 
Corollary \ref{4.5}.
\par\medskip

Because of the above results, we believe that approximation property 
is a good candidate for the notion of amenability
of actions of locally compact groups on general $C^*$-algebras. 
\par\medskip
\par\medskip

\noindent {\em II. Discrete groups: $G$-gradings and coactions.}
\par\medskip

Let $G$ be a discrete group and let $D$ be a $C^*$-algebra with a 
$G$-grading (i.e. $D=\overline{\oplus_{r\in G} D_r}$ such that 
$D_r\cdot D_s \subseteq D_{rs}$ and $D_r^*\subseteq D_{r^{-1}}$). 
Then there exists a canonical $C^*$-algebraic bundle structure 
(over $G$) on $D$. 
We denote this bundle by $\cal D$. 
Now by [6, \S 3], $D$ is a quotient of $C^*(\cal D)$.
Moreover, if the grading is topological in the sense that there exists a 
continuous conditional expectation from $D$ to $D_e$ (see [6, 3.4]), 
then $C^*_r({\cal D})$ is a quotient of $D$ (see [6, 3.3]). 
Hence by [12, 3.2(1)] (or [10, 2.17]), there is an induced 
non-degenerate coaction on $D$ by $C^*_r(G)$ 
which define the given grading. 
Now the proof of [9, 2.6], [6, 3.3] and the above 
observation imply the following equivalence. 
\par\medskip

\begin{proposition} \label{4.6} 
Let $G$ be a discrete group and $D$ be a $C^*$-algebra. 
Then a $G$-grading $D=\overline{\oplus_{r\in G} D_r}$ is topological 
if and only if it is induced by a non-degenerate coaction 
of $C^*_r(G)$ on $D$. 
\end{proposition}

\begin{corollary} \label{4.7}
Let $D$ be a $C^*$-algebra with a non-degenerate coaction $\epsilon$ 
by $C^*_r(G)$.
Then it can be ``lifted'' to a full coaction i.e. there exist a 
$C^*$-algebra $A$ with a full coaction $\epsilon_A$ by 
$G$ and a quotient map $q$ from $A$ to $D$ such that 
$\epsilon\circ q = (q\otimes \lambda_G)\circ\epsilon_A$. 
\end{corollary}

In fact, if $\cal D$ is the bundle as defined above, then we can take 
$A=C^*({\cal D})$ and $\epsilon_A = \delta_{\cal D}$.
\par\medskip
\par\medskip
\par\medskip

\noindent {\bf References}
\par\medskip
\noindent 
[1] F. Abadie, Tensor products of Fell bundles over discrete groups, 
preprint (funct-an/9712006), Universidade de S\~{a}o Paulo, 1997. 
\par\noindent
[2] C. Anantharaman-Delaroche, Action moyennable d'un groupe 
localement compact sur une alg\`ebre de von Neumann, 
Math. Scand. 45 (1979), 289-304. 
\par\noindent 
[3] C. Anantharaman-Delaroche, Action moyennable d'un groupe 
localement compact sur une alg\`ebre de von Neumann II, 
Math. Scand. 50 (1982), 251-268. 
\par\noindent 
[4] C. Anantharaman-Delaroche, Sys\`emes dynamiques non 
commutatifs et moyennabilit\`e, Math. Ann. 279 (1987), 297-315.
\par\noindent 
[5] S. Baaj and G. Skandalis, $C^{*}$-alg\`ebres de Hopf et 
th\'eorie de Kasparov \'equivariante, $K$-theory 2 (1989), 683-721.
\par\noindent 
[6] R. Exel, Amenability for Fell bundles, J. Reine Angew. Math., 
492 (1997), 41--73. 
\par\noindent 
[7] J. M. G. Fell and R. S. Doran, {\it Representations of *-algebras, 
locally compact groups, and Banach *-algebraic bundles vol. 1 and 2}, 
Academic Press, 1988.
\par\noindent 
[8] K. Jensen and K. Thomsen, {\it Elements of $KK$-Theory}, 
Birkh\"auser, 1991.
\par\noindent
[9] C. K. Ng, Discrete coactions on $C^*$-algebras, 
J. Austral. Math. Soc. (Series A), 60 (1996), 118-127.
\par\noindent 
[10] C. K. Ng, Coactions and crossed products of Hopf 
$C^{*}$-algebras, Proc. London. Math. Soc. (3), 72 (1996), 638-656.
\par\noindent 
[11] G. K. Pedersen, {\it $C^*$-algebras and their automorphism groups}, 
Academic Press, 1979.
\par\noindent 
[12] I. Raeburn, On crossed products by coactions and their 
representation theory, Proc. London Math. Soc. (3), 64 (1992), 625-652.
\par\noindent
[13] M. A. Rieffel, Induced representations of $C^*$-algebras, 
Adv. Math. 13 (1974), 176--257.
\par\noindent 
\par
\medskip
\noindent Departamento de Matem\'{a}tica, Universidade Federal de Santa
Catarina, 88010-970 Florian\'{o}polis SC, Brazil. 
\par
\noindent $E$-mail address: exel@mtm.ufsc.br
\par\medskip

\noindent Mathematical Institute, Oxford University, 24-29 St. Giles, Oxford 
OX1 3LB, United Kingdom.
\par
\noindent $E$-mail address: ng@maths.ox.ac.uk
\par
\end{document}